\documentclass[11pt,fleqn]{article}
\usepackage{amsmath,amssymb,graphicx}

\begin{document}
\sloppy
\newtheorem{axiom}{Axiom}[section]
\newtheorem{conjecture}[axiom]{Conjecture}
\newtheorem{corollary}[axiom]{Corollary}
\newtheorem{definition}[axiom]{Definition}
\newtheorem{example}[axiom]{Example}
\newtheorem{lemma}[axiom]{Lemma}
\newtheorem{observation}[axiom]{Observation}
\newtheorem{proposition}[axiom]{Proposition}
\newtheorem{theorem}[axiom]{Theorem}

\newcommand{\proof}{\emph{Proof.}\ \ }
\newcommand{\qed}{~~$\Box$}
\newcommand{\rz}{{\mathbb{R}}}
\newcommand{\nz}{{\mathbb{N}}}
\newcommand{\zz}{{\mathbb{Z}}}
\newcommand{\eps}{\varepsilon}
\newcommand{\cei}[1]{\lceil #1\rceil}
\newcommand{\flo}[1]{\left\lfloor #1\right\rfloor}
\newcommand{\seq}[1]{\langle #1\rangle}

%%%%%%%%%%%%%%%%%%%%%%%%%%%%%%%%%%%%%%%%%%%%%%
%%  GJW: macros
%%%%%%%%%%%%%%%%%%%%%%%%%%%%%%%%%%%%%%%%%%%%%%
\newcommand{\class}{\mbox{$\Delta_2$P}}
\newcommand{\bbb}{{\mathbb{B}}}

%%%%%%%%%%%%%%%%%%%%%%%%%%%%%%%%%%%%%%%%%%%%%%%%%%%%%%%%%%%%%%%%%%%%%%%%%%
%%%%%%%%%%%%%%%%%%%%%%%%%%%%%%%%%%%%%%%%%%%%%%%%%%%%%%%%%%%%%%%%%%%%%%%%%%
\title{{\bf Uniqueness in quadratic and hyperbolic\\0-1 programming problems}}
\author{
\sc Vladimir G.\ Deineko\thanks{{\tt Vladimir.Deineko@wbs.ac.uk}.
Warwick Business School, The University of Warwick, Coventry CV4 7AL, United Kingdom}
\and
\sc Bettina Klinz\thanks{{\tt klinz@opt.math.tu-graz.ac.at}.
Institut f\"ur Optimierung und Diskrete Mathematik, TU Graz, Steyrergasse 30, A-8010 Graz, Austria}
\and
\sc Gerhard J.\ Woeginger\thanks{{\tt gwoegi@win.tue.nl}.
Department of Mathematics and Computer Science, TU Eindhoven, P.O.\ Box 513,
5600 MB Eindhoven, Netherlands}
}
\date{}
\maketitle

\begin{abstract}
We analyze the question of deciding whether a quadratic or a hyperbolic 0-1 programming
instance has a unique optimal solution.
Both uniqueness questions are known to be NP-hard, but are unlikely to be contained in the class NP.
We precisely pinpoint their computational complexity by showing that they both are complete
for the complexity class {\class}.

\end{abstract}

%%%%%%%%%%%%%%%%%%%%%%%%%%%%%%%%%%%%%%%%%%%%%%%%%%%%%%%%%%%%%%%%%%%%%%%%%
%%%%%%%%%%%%%%%%%%%%%%%%%%%%%%%%%%%%%%%%%%%%%%%%%%%%%%%%%%%%%%%%%%%%%%%%%
\section{Introduction}
%%%%%%%%%%%%%%%%%%%%%%%%%%%%%%%%%%%%%%%%%%%%%%%%%%%%%%%%%%%%%%%%%%%%%%%%%
The general quadratic 0-1 programming problem has been studied by numerous authors; see for instance
Hammer \& Rudeanu \cite{HaRu1968} and Hansen \cite{Hansen1979}.
It has the form
%%%%%%%%%%%%%%%%
\begin{equation}
\label{eq:qua}
\min_{x\in\bbb^n}~~ x^TAx
\end{equation}
%%%%%%%%%%%%%%%%
where $A$ is an $n\times n$ symmetric integer matrix, and $\bbb^n$ is the set of all 0-1 vectors
$x=(x_1,x_2,\ldots,x_n)$ in dimension $n$.
As 0-1 variables satisfy $x_i^2=x_i$, the quadratic part of the objective function implicitly also
contains a linear function.
The quadratic 0-1 programming problem (\ref{eq:qua}) is well-known to be NP-hard, as -- for instance --
the maximum clique problem is a special case.
Hammer \& al \cite{Hammer2002} have shown that this problem remains NP-hard even if the objective function
in (\ref{eq:qua}) is the product of two linear 0-1 functions.

Another important special class of 0-1 optimization problems concerns the maximization of the ratio of
two linear functions with integer coefficients $a_0,\ldots,a_n$ and $b_0,\ldots,b_n$:
%%%%%%%%%%%%%%%%
\begin{equation}
\label{eq:hyp}
\max_{x\in\bbb^n}~~ \displaystyle \frac{a_0+\sum_{i=1}^na_ix_i}{b_0+\sum_{i=1}^nb_ix_i}.
\end{equation}
%%%%%%%%%%%%%%%%
This problem is also known as unconstrained single-ration hyperbolic 0-1 programming problem \cite{HaRu1968}.
If the denominator $b_0+\sum_{i=1}^nb_ix_i$ is positive for all $x\in\bbb^n$, then problem
(\ref{eq:hyp}) is polynomially solvable.
The problem in its full generality (where the denominator may take both negative and positive values)
however is NP-hard; see Boros \& Hammer \cite{BoHa2002}.

%%%%%%%%%%%%%%%%
\paragraph{Uniqueness of optimal solutions.}
Pardalos \& Jha \cite{Pardalos-1} discuss the complexity of deciding whether the quadratic optimization problem 
(\ref{eq:qua}) has a unique optimal solution $x^*\in\bbb^n$; they establish the NP-hardness of this question.
Prokopyev, Huang \& Pardalos \cite{Pardalos-2} analyze the complexity of deciding whether the hyperbolic
optimization problem (\ref{eq:hyp}) has a unique optimal solution, and prove also this question to be NP-hard.
(We note that the papers \cite{Pardalos-1,Pardalos-2} also show that quadratic $0$-$1$ and fractional $0$-$1$ 
programs remain NP-hard even if it is known that they have a unique optimal solution.)

There is no reason to believe that the problems of deciding uniqueness are contained in the complexity class NP 
(or in the closely related class coNP): to certify non-uniqueness of a concrete instance, one apparently 
has to exhibit two feasible solutions (that's the easy part which lies in NP) together with a certificate 
that these two solutions indeed are optimal (and that's the hard part).
Certifying that a solution is optimal amounts to proving the non-existence of a better solution,
which seems to require a coNP-certificate.
Such a mixture of NP- and coNP-certificates suggests that these questions might be located in one of the
complexity classes above NP and coNP; see for instance Chapter~{17} in Papadimitriou's book \cite{PapaBook}.

%%%%%%%%%%%%%%%%
\paragraph{Our results.}
We precisely pinpoint the computational complexity of the uniqueness questions for the quadratic program
(\ref{eq:qua}) and for the hyperbolic program (\ref{eq:hyp}): both problems are complete for the
complexity class {\class}.
The uniqueness questions for the quadratic program is {\class}-complete even if the objective function
is the product of two linear 0-1 functions.

The note is structured as follows.
Section~\ref{sec:prel} gives some preliminaries, and Section~\ref{sec:unique} presents our hardness proofs.

%%%%%%%%%%%%%%%%%%%%%%%%%%%%%%%%%%%%%%%%%%%%%%%%%%%%%%%%%%%%%%%%%%%%%%%%%
%%%%%%%%%%%%%%%%%%%%%%%%%%%%%%%%%%%%%%%%%%%%%%%%%%%%%%%%%%%%%%%%%%%%%%%%%
\section{Technical preliminaries}
\label{sec:prel}
%%%%%%%%%%%%%%%%%%%%%%%%%%%%%%%%%%%%%%%%%%%%%%%%%%%%%%%%%%%%%%%%%%%%%%%%%
An important and natural complexity class is {\class}, the class of all problems that can be solved
in polynomial time while using an oracle from NP; see \cite{PapaBook}.
It is known that {\class} contains the entire Boolean hierarchy over NP, that is, the smallest
complexity class that contains NP and that is closed under union, intersection, and complement.
Intuitively speaking (and assuming P$\ne$NP), this class {\class} is much larger than NP and
contains problems that are much more difficult than all the problems in NP and coNP.
A milestone result in this area by Papadimitriou \cite{Papadimitriou1984} shows that deciding whether a
given instance of the Travelling Salesman Problem (TSP) has a unique optimal solution is {\class}-complete;
in other words, the uniqueness question for the TSP belongs to the hardest problems in {\class} and
thus captures the full difficulty of the class {\class}.
As side results \cite{Papadimitriou1984} also established {\class}-completeness of deciding whether the 
optimal solution of an integer program or of a knapsack problem is unique; the side result on the knapsack
problem is the starting point of our proofs.

An instance to the (standard) knapsack problem consists of $n$ items with non-negative integer weights
$w_1,\ldots,w_n$ and non-negative integer profits $p_1,\ldots,p_n$, together with a weight bound $W$.
For a subset $I\subseteq\{1,2,\ldots,n\}$ we denote its weight by $w(I)=\sum_{i\in I}w_i$ and
its profit by $p(I)=\sum_{i\in I}p_i$; subset $I$ is called feasible if it satisfies $w(I)=W$.
By adding appropriate dummy items with profit~0, we may (and will) assume that there exists at least one
feasible subset.
The goal in the standard knapsack problem is to find a feasible subset $I$ that maximizes $p(I)$.

%%%%%%%%%%%%%%%
\begin{proposition}
\label{pr:Papadimitriou}
(Papadimitriou \cite{Papadimitriou1984})\\
It is {\class}-complete to decide whether a given instance of the knapsack problem has a
unique optimal solution.
\qed
\end{proposition}
%%%%%%%%%%%%%%%

In the proof of Theorem~\ref{th:qua}, we will work with the following variant of the knapsack problem 
which we call the \emph{subset sum problem with goal value}.
An instance consists of $n$ items with non-negative integer weights $q_1,\ldots,q_n$ together with
a goal value $Q$; as usual we denote $q(I)=\sum_{i\in I}q_i$ for an item set $I$.
The objective is to find an item subset $I$ that minimizes $|q(I)-Q|$ and hence has its weight
as close as possible to the goal value $Q$.

%%%%%%%%%%%%%%%
\begin{lemma}
\label{le:gjw}
It is {\class}-complete to decide whether a given instance of the subset sum problem with goal value
has a unique optimal solution.
\end{lemma}
%%%%%%%%%%%%%%%
\proof
Routine arguments show that the considered  uniqueness question is contained in {\class}; see for instance 
Lemma~\ref{le:containment} for similar arguments.
Hardness for {\class} is established by a reduction from the knapsack variant in Proposition~\ref{pr:Papadimitriou}.
Hence consider a knapsack instance, and let $P$ denote the total profit of all items.
For every item $i$ in the knapsack instance, we create a corresponding item $i$ in the subset sum instance with weight $q_i=3w_iP+p_i$.
Furthermore we define the goal value as $Q=(3W+1)P$.

Now consider an arbitrary item set $I$; note that $q(I)=3Pw(I)+p(I)$.
If $w(I)\le W-1$ then $q(I)\le Q-3P$, and if $w(I)\ge W+1$ then $q(I)\ge Q+2P$; either case yields a large objective value $|q(I)-Q|\ge2P$.
The remaining case $w(I)=W$ satisfies $Q-P\le q(I)\le Q$ and has a moderately small objective value $|q(I)-Q|\le P$.
Consequently the minimizers $I$ for the constructed subset sum instance satisfy $w(I)=W$ and hence are feasible for the knapsack instance.
Every such feasible set $I$ has $q(I)=Q-P+p(I)$, and thus minimizing $|q(I)-Q|=P-p(I)$ amounts to maximizing $p(I)$.
We conclude that the optimal solutions of the knapsack problem are in one-to-one correspondence with the optimal solutions of 
the subset sum problem.
This completes the hardness argument.
\qed

%%%%%%%%%%%%%%%%%%%%%%%%%%%%%%%%%%%%%%%%%%%%%%%%%%%%%%%%%%%%%%%%%%%%%%%%%
%%%%%%%%%%%%%%%%%%%%%%%%%%%%%%%%%%%%%%%%%%%%%%%%%%%%%%%%%%%%%%%%%%%%%%%%%
\section{Hardness of quadratic and hyperbolic programming}
\label{sec:unique}
%%%%%%%%%%%%%%%%%%%%%%%%%%%%%%%%%%%%%%%%%%%%%%%%%%%%%%%%%%%%%%%%%%%%%%%%%
In this section we analyze the uniqueness questions for the quadratic optimization problem (\ref{eq:qua})
and for the hyperbolic optimization problem (\ref{eq:hyp}).
We will first establish the containment of these questions in class {\class} and then provide the
{\class}-hardness proofs.

%%%%%%%%%%%%%%%
\begin{lemma}
\label{le:containment}
The uniqueness questions for the quadratic 0-1 programming problem and for the hyperbolic 0-1
programming problem are contained in class {\class}.
\end{lemma}
%%%%%%%%%%%%%%%
\proof
We first discuss the quadratic programming problem.
For an instance of (\ref{eq:qua}) let $S$ denote the sum of the absolute values of all entries in matrix $A$.
Note that $\log S$ is polynomially bounded in the instance size, and note that the optimal objective value
is an integer between $0$ and $S$.
In a first phase, we perform a binary search on the interval $[0,S]$ to determine the precise cost $c^*$
of the optimum by using an NP-oracle; this binary search takes at most $O(\log S)$ steps.
In the second phase we ask the NP-oracle whether there exist two vectors $x',x''\in\bbb^n$ that both yield
the objective value $c^*$.
If no such vectors exist we answer YES, and otherwise we answer NO.

The argument for the hyperbolic programming problem proceeds along the same lines.
For an instance of (\ref{eq:hyp}) let $S=\sum_{i=0}^n|a_i|+\sum_{i=0}^n|b_i|$.
Then $\log S$ is polynomially bounded in the instance size, and the optimal objective value is a rational
number whose numerator and denominator are integers between $0$ and $S$.
In a first phase, we determine the precise optimal cost $c^*$ by using an NP-oracle.
By applying standard methods from the literature, this search for a rational number can be implemented
in $O(\log S)$ steps; see Reiss \cite{Reiss1979} and in particular Zemel \cite{Zemel1981}.
The second phase again asks the NP-oracle for two vectors in $\bbb^n$ that both have objective value $c^*$.
\qed

\bigskip
Next, we will establish {\class}-hardness of the uniqueness question for the quadratic 0-1 programming problem.
In fact we will prove {\class}-hardness even for the highly restricted special case where the objective function 
is the product of two linear 0-1 functions so that matrix $A$ in (\ref{eq:qua}) has rank~1.
The proof uses a polynomial time reduction from the subset sum problem in Lemma~\ref{le:gjw}.
Hence consider a subset sum instance with weights $q_1,\ldots,q_n$ and goal value $Q$. 
We define the quadratic function 
%%%%%%%%%%%%%%%%
\begin{equation}
\label{eq:q1}
f(x) ~=~ \left( \sum_{i=1}^nq_ix_i -3Qx_{n+1}\right)\, \left(\sum_{i=1}^nq_ix_i +Qx_{n+1}\right).
\end{equation}
%%%%%%%%%%%%%%%%
Clearly, $f(x)$ is the product of two linear 0-1 functions.
Setting $x_i=0$ for $i=1,\ldots,n$ and setting $x_{n+1}=1$ yields the function value $f(x)=-3Q^2<0$; 
hence the optimal objective value will be negative.
For $x_{n+1}=0$, the function $f(x)$ is a perfect square and only takes non-negative values;
hence these cases can never yield a minimizer of $f$.
For $x_{n+1}=1$, the function $f(x)$ can be rewritten as
%%%%%%%%%%%%%%%%
\begin{equation}
\label{eq:q2}
f(x) ~=~ \left( \sum_{i=1}^nq_ix_i -3Q\right)\, \left(\sum_{i=1}^nq_ix_i +Q\right)
~=~ \left( \sum_{i=1}^nq_ix_i -Q\right)^2-4Q^2.
\end{equation}
%%%%%%%%%%%%%%%%
Hence, in this case the quadratic function $f(x)$ has exactly the same minimizers as the function 
$\left|\sum_{i=1}^nq_ix_i -Q\right|$ in the subset sum instance.
This one-to-one correspondence yields the following theorem.

%%%%%%%%%%%%%%%
\begin{theorem}
\label{th:qua}
It is {\class}-complete to decide whether a given instance of the quadratic 0-1 programming problem
has a unique optimal solution.
This hardness result holds even for the special case where the objective function is the product of 
two linear 0-1 functions (and where matrix $A$ is of rank~1).
\qed
\end{theorem}
%%%%%%%%%%%%%%%

\bigskip
Now let us show {\class}-hardness of the uniqueness question for the hyperbolic 0-1 programming problem
through a reduction from the knapsack variant in Proposition~\ref{pr:Papadimitriou}.
Consider a knapsack instance, and let $P$ denote the total profit of all items.
For a vector $x=(x_1,x_2,\ldots,x_n)$ in $\bbb^n$ we consider the hyperbolic function
%%%%%%%%%%%%%%%%
\begin{equation}
\label{eq:h}
f(x) ~=~  \left( 2+\sum_{i=1}^np_ix_i \right) ~\mbox{\LARGE/}~ \left( 1+2P\cdot (W-\sum_{i=1}^nw_ix_i)\right)
\end{equation}
%%%%%%%%%%%%%%%%
The knapsack instance has at least one feasible item set, and the corresponding characteristic vector
$x\in\bbb^n$ satisfies $\sum_{i=1}^nw_ix_i=W$.
For this vector $x$, the denominator in (\ref{eq:h}) evaluates to~$1$; hence $f(x)\ge2$ and the
maximum objective value of (\ref{eq:h}) is at least $2$.
If some vector $x$ satisfies $\sum_{i=1}^nw_ix_i<W$, then the denominator in (\ref{eq:h}) is at least $1+2P$,
whereas the numerator is at most $2+P$; hence $f(x)<2$.
If some vector $x$ satisfies $\sum_{i=1}^nw_ix_i>W$, then the denominator in (\ref{eq:h}) is negative
and the numerator is positive; hence $f(x)<0$.
Altogether this implies that every maximizing solution of (\ref{eq:h}) will satisfy $\sum_{i=1}^nw_ix_i=W$.

We conclude that the maximizers of (\ref{eq:h}) are the maximizers of the numerator $2+\sum_{i=1}^np_ix_i$.
Consequently, a vector $x$ maximizes (\ref{eq:h}) if and only if the corresponding item set $I$ is feasible
and maximizes the profit.
Hence, optimal solutions of the knapsack problem are in one-to-one correspondence with optimal solutions
of the hyperbolic problem (\ref{eq:h}), and we get the following theorem.

%%%%%%%%%%%%%%%
\begin{theorem}
\label{th:hyp}
It is {\class}-complete to decide whether a given instance of the hyperbolic 0-1 programming problem
has a unique optimal solution.
\qed
\end{theorem}
%%%%%%%%%%%%%%%

Finally, let us spend a few words on the consequences of our result for 0-1 programming.
If some problem $X$ is {\class}-complete, then it encapsulates the full difficulty of class {\class}.
This means that problem $X$ is essentially equally powerful as an entire sequence of (polynomially many)
integer programs, where the formulations of the later IPs in the sequence may depend on the solutions of
the earlier IPs in the sequence.
This is of course strong evidence that problem $X$ cannot be expressed as a single standard IP.
Therefore,  any {\class}-complete problem $X$ (and in particular the uniqueness questions for quadratic 
and hyperbolic 0-1 programs) is well beyond reach of the classical mathematical programming machinery.

%%%%%%%%%%%%%%%%%%%%%%%%%%%%%%%%%%%%%%%%%%%%%%%%%%%%%%%%%%%%%%%%%%%%%%%%%
%%%%%%%%%%%%%%%%%%%%%%%%%%%%%%%%%%%%%%%%%%%%%%%%%%%%%%%%%%%%%%%%%%%%%%%%%
\section*{Acknowledgements}

Vladimir Deineko acknowledges support
by Warwick University's Centre for Discrete Mathematics and Its Applications (DIMAP)
and by EPSRC fund EP/F017871.
Bettina Klinz acknowledges partial financial support by the
FWF program W 1230-N13.
Gerhard Woeginger acknowledges support
by the Netherlands Organization for Scientific Research (NWO), grant 639.033.403,
and by DIAMANT (an NWO mathematics cluster).

%%%%%%%%%%%%%%%%%%%%%%%%%%%%%%%%%%%%%%%%%%%%%%%%%%%%%%%%%%%%%%%%%%%%%%%%%
%%%%%%%%%%%%%%%%%%%%%%%%%%%%%%%%%%%%%%%%%%%%%%%%%%%%%%%%%%%%%%%%%%%%%%%%%

\end{document}